\documentclass[10pt]{article}

\usepackage{amsthm}
\usepackage{amsmath}
\usepackage{amsfonts}
\usepackage{amssymb}
\usepackage{mathrsfs}
\usepackage{empheq}
\usepackage{stmaryrd}
\usepackage{enumerate}
\usepackage{pgfplots,soul}
\usepackage[applemac]{inputenc}
\usepackage[colorlinks=true,urlcolor=blue]{hyperref}

\newcommand{\N}{\mathbb{N}}

\newcommand{\R}{\mathbb{R}}
\newcommand{\C}{\mathbb{C}}

\newcommand{\Dr}{\mathscr{D}}
\newcommand{\Fr}{\mathscr{F}}

\newcommand{\vphi}{\varphi}
\newcommand{\eps}{\varepsilon}

\newcommand{\dsp}{\displaystyle}
\newcommand{\ovl}{\overline}
\newcommand{\udl}{\underline}
\newcommand{\vlim}{\lim\limits}

\newcommand{\vint}{\int\limits}

\newcommand{\inj}{\hookrightarrow}
\newcommand{\tends}{\longrightarrow}
\newcommand{\weak}{\rightharpoonup}

\newcommand{\loc}{\mathrm{loc}}

\renewcommand{\b}{\mathrm{b}}

\renewcommand{\d}{\mathrm{d}}

\newcommand{\vi}{\mathrm{i}}
\newcommand{\w}{{\textsl w}}
\renewcommand{\le}{\leqslant}
\renewcommand{\ge}{\geqslant}
\renewcommand{\Re}{\mathrm{Re}}
\renewcommand{\Im}{\mathrm{Im}}
\newcommand{\bs}{\boldsymbol}
\newcommand{\p}{\prime}

\newcommand{\eqdef}{\stackrel{\mathrm{def}}{=}}
\DeclareMathOperator{\supp}{supp}

\DeclareMathOperator{\Arg}{Arg}

\numberwithin{equation}{section}

\newtheorem{thm}{Theorem}[section]

\newtheorem{cor}[thm]{Corollary}
\newtheorem{lem}[thm]{Lemma}

\theoremstyle{definition}
\newtheorem{rmk}[thm]{Remark}
\newtheorem{defi}[thm]{Definition}

\newtheorem{assum}[thm]{Assumption}

\newenvironment{proof*}{\noindent{\bf Proof.}}{\qed}
\newenvironment{vproof}[1]{\noindent{\bf Proof #1}}{\qed}

\title{\huge \sc Finite time extinction for a damped nonlinear Schrödinger equation in the whole space}
\author{\sc Pascal Bégout}
\date{}

\begin{document}

\maketitle

\begin{center}
Institut de Mathématiques de Toulouse \& TSE	\\
Université Toulouse I Capitole				\\
1, Esplanade de l’Université				\\
31080 Toulouse Cedex 6, FRANCE
\bigskip \\
\text{
{\footnotesize E-mail\:: }\htmladdnormallink{{\footnotesize\udl{\tt{Pascal.Begout@math.cnrs.fr}}}}{mailto:Pascal.Begout@math.cnrs.fr}}
\end{center}

\begin{abstract}
We consider a nonlinear Schrödinger equation set in the whole space with a single power of interaction and an external source. We first establish existence and uniqueness of the solutions and then show, in low space dimension, that the solutions vanish at a finite time. Under a smallness hypothesis of the initial data and some suitable additional assumptions on the external source, we also show that we can choose the upper bound on which time the solutions vanish.
\end{abstract}

{\let\thefootnote\relax}
{\let\thefootnote\relax}

\tableofcontents

\baselineskip .6cm

\section{Introduction and explanation of the method}
\label{introduction}

Let us consider the following Schrödinger equation with a nonlinear damping term,
\begin{gather}
\label{nlsi}
\vi u_t+\Delta u+a|u|^{m-1}u=f(t,x), \; \text{ in } (0,\infty)\times\Omega,
\end{gather}
where $\Omega\subseteq\R^N$ is an open subset, $a\in\C,$ $0<m<1$ and $f:(0,\infty)\times\Omega\tends\C$ measurable is an external source. When $a\in\R,$ $m\ge1$ and $f=0,$ equation~\eqref{nlsi} has been intensively studied, especially with $\Omega=\R^N$ (among which existence, uniqueness, blow-up, scattering theory, time decay). The literature is too extensive to give an exhaustive list. See, for instance, the monographs of Cazenave~\cite{MR2002047}, Sulem and Sulem~\cite{MR2000f:35139}, Tao~\cite{MR2233925} and the references therein. The case $a\in\C$ is more anecdotic. See, for instance, Bardos and Brezis~\cite{MR0242020}, Lions~\cite{MR0259693}, Tsutsumi ~\cite{MR1038160} and Shimomura~\cite{MR2254620}. Note that except in \cite{MR0259693}, it is always assumed $m>1.$
\medskip
\\
In this paper, we are looking for solutions which vanishes at a finite time. For many reasons, we have to consider $0<m<1.$ When $m=1,$ existence is not hard to obtain, since the equation is linear, while the finite time property is not possible (which is a direct consequence of~\eqref{introM}). To our knowledge the first paper in this direction is due to Carles and Gallo~\cite{MR2765425} with $a=\vi,$ $f=0$ and $\Omega$ is a compact manifold without boundary. To construct solutions, they regularize the nonlinearity and use a compactness method to pass in the limit. They prove the finite time extinction property for $N\le3$ including the case $m=0.$ More recently, Carles and Ozawa~\cite{MR3306821} obtain the existence, uniqueness and finite time extinction for $\Omega=\R^N,$ $a\in\vi\R_+$ and $f=0.$ Due to the lack of compactness, they restrict their study to $N\le2$ and add an harmonic confinement in~\eqref{nlsi} for some technical reasons. For the finite time property with $N=2$ they also restrict the range of $m$ to $\left[\frac12,1\right)$ and make a smallness assumption of the initial data. In this paper, we work in the whole space and we remove of all these restrictions and extend the previous results to a large class of values of $a$ (see, for instance, Theorems~\ref{thmstrongH2} and \ref{thmextH2}). Indeed, we shall assume that the complex number $a$ is in a cone of the complex plane. More precisely,
\begin{gather}
\label{introa}
a\in C(m)\eqdef\Big\{z\in\C; \; \Im(z)>0 \; \text{ and } \; 2\sqrt m\Im(z)\ge(1-m)|\Re(z)|\Big\}.
\end{gather}
The assumption that $a$ belongs to the cone $C(m)$ was considered in a series of papers by Okazawa and Yokota~\cite{MR1970632,MR1900334,MR1886827}. They studied the asymptotic behavior of the solutions to the complex Ginzburg-Landau equation in a bounded domain with the assumption~\eqref{introa} and, sometimes, with $m>1.$ See also Kita and Shimomura~\cite{MR2272871} and Hou, Jiang, Li and You~\cite{MR3228869} where~\eqref{introa} is assumed but with (among others restrictive assumptions) $m>1.$ In all these papers, there is no finite time extinction result. We would also like mention the (very complete) work of Antontsev, Dias and Figueira~\cite{MR3208711} where they consider the complex Ginzburg-Landau equation,
\begin{gather}
\label{gl}
e^{-\vi\gamma}u_t-\Delta u+|u|^{m-1}u=f(t,x), \; \text{ in } (0,\infty)\times\Omega,
\end{gather}
where $\Omega$ is bounded, $0<m<1$ and $-\frac\pi2<\gamma<\frac\pi2.$ In particular, $e^{-\vi\gamma}\neq\pm\vi.$ They show spatial localization, waiting time and finite time extinction properties. The case of equation~\eqref{gl} with a delayed nonlocal perturbation is studied in the recent paper of D\'iaz, Padial, Tello and Tello~\cite{dptt}. Finally, Hayashi, Li and Naumkin~\cite{MR3465033} study time decay for a more classical Schrödinger equation~\eqref{nlsi} $(a$ satisfying~\eqref{introa}, $m>1$ and $\Omega=\R^N).$
\medskip
\\
In this paper, we are interested in the finite time extinction of the solution. Formally, this result is not too hard to obtain (the method we explain below for the finite time extinction property is that used in \cite{MR2765425,MR3306821,MR4053613}). Suppose $f=0.$ It is well known that solutions that vanish in finite time do not exist when $m\ge1$ (at least when $a\in\R).$ Indeed, multiplying~\eqref{nlsi} by $\ovl{\vi u},$ integrating by parts and taking the real part, we obtain,
\begin{gather}
\label{introM}
\frac12\frac{\d}{\d t}\|u(t)\|_{L^2}^2+\Im(a)\|u(t)\|_{L^{m+1}}^{m+1}=0.
\end{gather}
To expect a finite time extinction, the mass has to be non increasing and so $\Im(a)>0.$ Now, since $m+1<2,$ we may interpolate $L^2$ between $L^{m+1}$ and $L^p,$ for some $p>2,$ and control the $L^p$-norm by a Sobolev norm. Using a Gagliardo-Nirenberg's inequality,
\begin{gather}
\label{intro1}
\|u(t)\|_{L^2}^{2\frac{m+1}{2\theta_\ell}}\le\|u(t)\|_{L^{m+1}}^{m+1}\|u(t)\|_{H^\ell}^\frac{(m+1)(1-\theta_\ell)}{\theta_\ell},
\end{gather}
for some an explicit constant $\theta_\ell\in(0,1),$ if $u$ is bounded in $H^\ell$ then putting together \eqref{introM}--\eqref{intro1}, we arrive at the ordinary differential equation,
\begin{gather}
\label{intro2}
y^\p+Cy^\delta\le0,
\end{gather}
with $\delta=\frac{m+1}{2\theta_\ell},$ where $y(t)=\|u(t)\|_{L^2}^2.$ By integration, we then obtain the asymptotic behavior of $u$ with respect to the value of $\delta.$

$\bullet$
If $\delta<1$ then $y(t)^{1-\delta}\le(y(0)^{1-\delta}-Ct)_+$ and so $u$ vanishes before time $T_\star=C^{-1}y(0)^{1-\delta}.$

$\bullet$
If $\delta=1$ then $y(t)\le y(0)e^{-Ct}.$

$\bullet$
If $\delta>1$ then $y(t)^{\delta-1}\le y(0)^{\delta-1}(1+Ct)^{-1}.$

\noindent
As a consequence, a sufficient condition to have extinction in finite time is $\delta<1$ which turns out to be equivalent to $N=1$ when $\ell=1.$ To increase the space dimension, we assume that $u$ is bounded in $H^2$ and we deduce that $\delta<1$ when $N\le3.$ Theoretically, we can reach any space dimension if $u$ is bounded in $H^\ell$ for $\ell$ large enough (actually, if $\ell=\left[\frac{N}2\right]+1,$ where $\left[\frac{N}2\right]$ denotes the integer part of $\frac{N}2;$ see Theorem~2.1 in Bégout and D\'iaz~\cite{MR4053613}). But this is not reasonable due to the lack of regularity of the nonlinearity, which is merely Hölder continuous. A reachable goal is to obtain existence and boundedness of the solutions in $H^2.$
\medskip
\\
Now, we focus on the construction of a solution to~\eqref{nlsi} in $\R^N$ with $f=0$ (to fix ideas). First of all, we would like to uniformly control $\|u(t)\|_{H^1}^2.$ Estimate~\eqref{introM} partially answers this question. For $\|\nabla u(t)\|_{L^2}^2,$ we multiply~\eqref{nlsi} by $\vi\ovl{\Delta u}$ and take the real part. We get,
\begin{gather*}
\frac12\frac{\d}{\d t}\|\nabla u(t)\|_{L^2}^2+\Re\left(\vi a\vint_{\R^N}|u(t)|^{m-1}u(t)\ovl{\Delta u(t)}\d x\right)=0.
\end{gather*}
We then expect to have,
\begin{gather}
\label{introE}
\Re\left(\vi a\vint_{\R^N}|u(t)|^{m-1}u(t)\ovl{\Delta u(t)}\d x\right)\ge0.
\end{gather}
Regularizing the nonlinearity, integrating by parts and passing to the limit,~\eqref{introE} can be proved under assumption~\eqref{introa} (Lemma~\ref{lemLap}). Actually, we extended the method found in Carles and Gallo~\cite{MR2765425}, where the situation is simpler since $a=\vi.$ Assume $\Omega\subseteq\R^N.$ To construct a solution to~\eqref{nlsi}, we use theory of the maximal monotone operators in the Hilbert space $L^2.$ We then consider the operator,
\begin{gather}
\label{introA}
Au=-\vi\Delta u-\vi a|u|^{m-1}u,
\end{gather}
with the natural domain\footnote{\label{fn1}It is natural in the sense that it is the smallest domain, in the sense of the inclusion, for which $D(A)\subset L^2.$} $D(A)=\big\{u\in H^1_0(\Omega);u^m\in L^2(\Omega)\text{ and } \Delta u\in L^2(\Omega)\big\}.$ Monotonicity relies on the inequality,
\begin{gather}
\label{introu}
\Re\left(-\vi\,a\vint_\Omega\big(|u|^{m-1}u-|v|^{m-1}v\big)(\ovl{u-v})\d x\right)\ge0.
\end{gather}
Once~\eqref{introu} is proved, it remains to show that $R(I+A)=L^2$ (Theorem~\ref{thmexicom} and Corollary~\ref{Amaxmon}). This means that for any $F\in L^2,$ the equation
\begin{gather}
\label{introF}
-\vi\Delta u-\vi a|u|^{m-1}u+u=F,
\end{gather}
admits a solution belonging to $D(A).$ Existence, uniqueness, \textit{a priori} estimates and smoothness of the solutions of~\eqref{introF} for a large class of values of $a$ (including \eqref{introa}) have been intensively studied in the papers by Bégout and D\'iaz~\cite{MR2876246,MR3315701}. The natural\footnote{\label{fn2}Multiply \eqref{introF} by $\ovl{\vi u}$ and $\ovl u,$ integrate by parts and take the real part.} space to look for a solution is $H^1_0\cap L^{m+1}.$ When $\Omega$ is bounded with a smooth boundary, a bootstrap method yields $u\in H^2(\Omega).$ Note that in this case, the condition $u^m\in L^2(\Omega)$ is automatically verified since $u^m\in L^\frac2m(\Omega)\inj L^2(\Omega)$ and then $u\in D(A).$ Although this method works very well, we proposed another one in Bégout and D\'iaz~\cite{MR4053613}: we make the sum of two monotone operators, where one of them is maximal monotone $(-\vi\Delta)$ and the other one is continuous over $L^2(\Omega)$ $(-\vi a|u|^{m-1}u).$ A difficulty appears when $\Omega$ is unbounded, say $\Omega=\R^N.$ In this case, we have $D(A)=H^2(\R^N)\cap L^{2m}(\R^N)$ and we have to show that a solution $u\in H^1(\R^N)\cap L^{m+1}(\R^N)$ belongs to $L^{2m}(\R^N),$ or equivalently $\Delta u\in L^2(\R^N).$ Having~\eqref{introE} in mind, a natural method would be to multiply~\eqref{introF} by $-\ovl{\Delta u}$ and take the real part. But then we lose the term $\|\Delta u\|_{L^2(\R^N)}^2.$ The original idea is to rotate $a$ in the complex plane and stay in the cone $C(m)$ to still have~\eqref{introE} (see Lemma~\ref{lemb} and the picture p.\pageref{picture}). If we can find $b\in\C$ such that $ab\in C(m)$ then multiplying~\eqref{introF} by $-b\ovl{\Delta u},$ integrating by parts and taking the real part, we arrive at,
\begin{gather*}
-\Im(b)\|\Delta u\|_{L^2(\R^N)}^2+\Re\left(\vi ab\vint_{\R^N}|u|^{m-1}u\ovl{\Delta u}\d x\right)+\Re(b)\|\nabla u\|_{L^2(\R^N)}^2
=-\Re\left(b\vint_{\R^N}F\ovl{\Delta u}\d x\right).
\end{gather*}
We see that we must have $\Im(b)<0$ and so the rotation has to be made in the negative sense. So we exclude the boundary of $C(m)$ located in the first quarter complex plane. Hence Assumption~\ref{ass} below. Note that the sign of $\Re(b)$ has no importance since we already have an estimate in $H^1(\R^N).$ Having \textit{a priori} estimates, we may construct a solution $u\in H^2(\R^N)\cap L^{2m}(\R^N)$ of~\eqref{introF} as a limit of solutions with compact support. The existence of such solutions is provided in Bégout and D\'iaz~\cite{MR2876246} (see also Bégout and D\'iaz~\cite{MR3190983}). To conclude the explanation of our method, we go back to the proof of~\eqref{introu}. When $a=\vi,$ this is very simple since this estimate is equivalent to the monotonicity of the derivative of the convex function defined on $\R^2$ by, $(x,y)\longmapsto\frac1{m+1}(x^2+y^2)^\frac{m+1}2$ (see Remark~9.3 in Bégout and D\'iaz~\cite{MR2876246}). But when $\Re(a)\neq0$ then the imaginary part of the integral in~\eqref{introu} is still there. Fortunately, this can be controlled by its real part under assumption~\eqref{introa} and a consequence of Liskevich and Perel$^\p$muter~\cite{MR1224619} (Lemma~2.2).
\medskip
\\
Finally, we consider the limit cases $m=0$ and $m=1$ for the values of $a.$ Since $\vlim_{m\searrow0}C(m)=\{0\}\times\vi(0,\infty),$ it seems that no extension of~\cite{MR2765425,MR3306821} is possible. The other limit case $\vlim_{m\nearrow1}C(m)=\R\times\vi(0,\infty)$ is entirely treated in Bégout and D\'iaz~\cite{MR4053613}: existence, uniqueness and boundedness for any subset $\Omega\subseteq\R^N.$
\medskip
\\
We will use the following notations throughout this paper. We denote by $\ovl z$ the conjugate of the complex number $z,$ by $\Re(z)$ its real part and by $\Im(z)$ its imaginary part. Unless if specified, all functions are complex-valued $(H^1(\Omega)=H^1(\Omega;\C),$ etc). For $1\le p\le\infty,$ $p^\prime$ is the conjugate of $p$ defined by $\frac{1}{p}+\frac{1}{p^\prime}=1.$ For a Banach space $X,$ we denote by $X^\star$ its topological dual and by $\langle\: . \; , \: . \:\rangle_{X^\star,X}\in\R$ the $X^\star-X$ duality product. In particular, for any $T\in L^{p^\prime}(\Omega)$ and $\vphi\in L^p(\Omega)$ with $1\le p<\infty,$ $\langle T,\vphi\rangle_{L^{p^\prime}(\Omega),L^p(\Omega)}=\Re\int_{\Omega}T(x)\ovl{\vphi(x)}\d x.$ The scalar product in $L^2(\Omega)$ between two functions $u,v$ is, $(u,v)_{L^2(\Omega)}=\Re\int_{\Omega}u(x)\ovl{v(x)}\d x.$ For a Banach space $X$ and $p\in[1,\infty],$ $u\in L^p_\loc\big([0,\infty);X\big)$ means that for any $T>0,$ $u_{|(0,T)}\in L^p\big((0,T);X\big).$ In the same way, we will use the notation $u\in W^{1,p}_\loc\big([0,\infty);X\big).$ As usual, we denote by $C$ auxiliary positive constants, and sometimes, for positive parameters $a_1,\ldots,a_n,$ write as $C(a_1,\ldots,a_n)$ to indicate that the constant $C$ depends only on $a_1,\ldots,a_n$ and that dependence is continuous (we will use this convention for constants which are not denoted by ``$C$'').
\medskip
\\
This paper is organized as follows. In Section~\ref{exiuni}, we state the mains results about existence, uniqueness and boundness for~\eqref{nlsi} (Theorem~\ref{thmweak}, \ref{thmstrongH1} and \ref{thmstrongH2}). In Section~\ref{finite}, we give the results about the finite time extinction property and the asymptotic behavior (Theorems~\ref{thmextH2}, \ref{thm0s} and \ref{thm0w}). The proofs of the existence, uniqueness and boundness are made in Section~\ref{proofexi} while those of the finite time extinction property and the asymptotic behavior are given in Section~\ref{proofext}.

\section{Existence and uniqueness of the solutions}
\label{exiuni}

Let $0<m<1,$ let $a\in\C,$ let $f\in L^1_\loc\big([0,\infty);L^2(\R^N)\big)$ and let $u_0\in L^2(\R^N).$ We consider the following nonlinear Schrödinger equation.

\begin{empheq}[left=\empheqlbrace]{align}
	\label{nls}
	\vi\frac{\partial u}{\partial t}+\Delta u+a|u|^{-(1-m)}u=f(t,x),	&	\text{ in } (0,\infty)\times\R^N,		\\
	\label{u0}
	u(0)= u_0,										&	\text{ in } \R^N,
\end{empheq}

\bigskip
\noindent
The main results in this paper hold with the assumptions below.

\begin{assum}
\label{ass}
We assume that $0<m<1$ and $a\in\C$ satisfy,
\begin{gather}
\label{a}
2\sqrt m\,\Im(a)\ge(1-m)|\Re(a)|.
\end{gather}
If $\Re(a)\ge0$ then we assume further that,
\begin{gather}
\label{as}
2\sqrt m\,\Im(a)>(1-m)\Re(a).
\end{gather}
\end{assum}
\medskip

\noindent
Here and after, we shall always identify $L^2(\R^N)$ with its topological dual. Let $0<m<1$ and let $X=H\cap L^{m+1}(\R^N),$ where $H=L^2(\R^N)$ or $H=H^1(\R^N).$ We recall that (see, for instance, Lemmas~A.2 and A.4 in Bégout and D\'iaz~\cite{MR4053613}),
\begin{gather}
\label{defX*}
X^\star=H^\star+L^\frac{m+1}m(\R^N),										\\
\label{X}
\Dr(\R^N)\inj X\inj L^{m+1}(\R^N) \text{ with both dense embeddings,}				\\
\label{X*}
L^\frac{m+1}m(\R^N)\inj X^\star\inj\Dr^\p(\R^N), \text{ with both dense embeddings,}		\\
\label{Xcon}
L^{m+1}_\loc\big([0,\infty);X\big)\cap W^{1,\frac{m+1}{m}}_\loc\big([0,\infty);X^\star\big)\inj C\big([0,\infty);L^2(\R^N)\big).
\end{gather}
This justifies the notion of solution below (and especially~\ref{defsol4})).

\begin{defi}
\label{defsol}
Let $0<m<1,$ let $a\in\C,$ let $f\in L^1_\loc\big([0,\infty);L^2(\R^N)\big)$ and let $u_0\in L^2(\R^N).$ Let us consider the following assertions.
\begin{enumerate}[1)]
\item
\label{defsol1}
$u\in L^{m+1}_\loc\big([0,\infty);H^1(\R^N)\cap L^{m+1}(\R^N)\big)
\cap W^{1,\frac{m+1}{m}}_\loc\big([0,\infty);H^\star+L^\frac{m+1}m(\R^N)\big),$
\item
\label{defsol2}
For almost every $t>0,$ $\Delta u(t)\in H^\star.$
\item
\label{defsol3}
$u$ satisfies~\eqref{nls} in $\Dr^\p\big((0,\infty)\times\R^N\big).$
\item
\label{defsol4}
$u(0)=u_0.$
\end{enumerate}
We shall say that $u$ is a \textit{strong solution} if $u$ is an $H^2$-solution or an $H^1$-solution. We shall say that $u$ is an $H^2$-\textit{solution of} \eqref{nls}--\eqref{u0} \big(respectively, an $H^1$-\textit{solution of} \eqref{nls}--\eqref{u0}\big), if $u$ satisfies the Assertions~\ref{defsol1})--\ref{defsol4}) with $H=L^2(\R^N)$ \big(respectively, with $H=H^1(\R^N)\big).$
\\
We shall say that $u$ is an $L^2$-\textit{solution} or a \textit{weak solution} of \eqref{nls}--\eqref{u0} is there exists a pair,
\begin{gather}
\label{fn}
(f_n,u_n)_{n\in\N}\subset L^1_\loc\big([0,\infty);L^2(\R^N)\big)\times C\big([0,\infty);L^2(\R^N)\big),
\end{gather}
such that for any $n\in\N,$ $u_n$ is an $H^2$-solution of \eqref{nls} where the right-hand side of \eqref{nls} is $f_n,$ and if
\begin{gather}
\label{cv}
f_n\xrightarrow[n\to\infty]{L^1((0,T);L^2(\R^N))}f \; \text{ and } \; u_n\xrightarrow[n\to\infty]{C([0,T];L^2(\R^N))}u,
\end{gather}
for any $T>0,$ and if $u$ satisfies~\eqref{u0}.
\end{defi}

\begin{rmk}
\label{rmkdefsol}
Let $0<m<1.$ Set for any $z\in\C,$ $g(z)=|z|^{-(1-m)}z$ $(g(0)=0).$ We define the mapping for any measurable function $u:\R^N\tends\C,$ which we still denote by $g,$ by $g(u)(x)=g(u(x)).$ Let $X$ be as in the beginning of this section (see~\eqref{defX*}--\eqref{Xcon}). From~\eqref{X}, \eqref{X*} and the basic estimate,
\begin{gather}
\label{g}
\forall(z_1,z_2)\in\C^2, \; |g(z_1)-g(z_2)|\le C|z_1-z_2|^m,
\end{gather}
(see, for instance, Lemma~A.1 in Bégout and D\'iaz~\cite{MR4053613}), we deduce easily that,
\begin{gather}
\label{gmmb}
g\in C\big(L^{m+1}(\R^N);L^\frac{m+1}m(\R^N)\big) \text{ and } g \text{ is bounded on bounded sets,}	\\
\label{gmm}
g\in C(X;X^\star) \text{ and } g \text{ is bounded on bounded sets.}
\end{gather}
By~\eqref{X}--\eqref{X*} and~\eqref{gmmb}--\eqref{gmm}, it follows that,
\begin{gather}
\label{dualg}
\langle g(u),v\rangle_{X^\star,X}=\langle g(u),v\rangle_{L^\frac{m+1}m(\R^N),L^{m+1}(\R^N)}
=\Re\vint_{\R^N} g(u)\ovl v\d x,
\end{gather}
for any $u,v\in X.$ Now, let us collect some basic informations about the solutions.
\begin{enumerate}[1)]
\item
\label{rmkdefsol1}
Any strong or weak solution belongs to $C\big([0,\infty);L^2(\R^N)\big)$ and Assertion~\ref{defsol4}) makes sense in $L^2(\R^N)$ (by~\eqref{Xcon}). 
\item
\label{rmkdefsol2}
It is obvious that an $H^2$-solution is also an $H^1$-solution and a weak solution. But it is not clear that an $H^1$-solution is a weak solution, without a continuous dependence of the solution with respect to the initial data. Such a result will be established with the additional assumptions~\eqref{a}--\eqref{as} on $a$ (see Lemma~\ref{lemdep} below). Note also that Assertion~\ref{defsol2}) of Definition~\ref{defsol} is not an additional assumption for the $H^1$-solutions.
\item
\label{rmkdefsol3}
Any $H^2$-solution (respectively, any $H^1$-solution) satisfies~\eqref{nls} in $L^2(\R^N)+L^\frac{m+1}m(\R^N)$ \big(respectively, in $H^{-1}(\R^N)+L^\frac{m+1}m(\R^N)\big),$ for almost every $t>0.$ Indeed, this is a direct consequence of Definition~\ref{defsol} and~\eqref{gmm}.
\item
\label{rmkdefsol4}
If $u$ is a weak solution then $u\in W^{1,1}_\loc\big([0,\infty);Y^\star\big)$ and it solves~\eqref{nls} in $Y^\star,$ for almost every $t>0,$ where $Y=H^2(\R^N)\cap L^\frac2{2-m}(\R^N)$ and $Y^\star=H^{-2}(\R^N)+L^\frac2m(\R^N)\inj\Dr^\p(\R^N)$ (by Lemma~A.2 in Bégout and D\'iaz~\cite{MR4053613}). Indeed, using the notation of Definition~\ref{defsol} and \eqref{g}, this comes from \eqref{cv} and the uniform convergences,
\begin{align}
\label{rmkdefsol41}
&	\Delta u_n\xrightarrow[n\to\infty]{C([0,T];H^{-2}(\R^N))}\Delta u,					\\
\label{rmkdefsol42}
&	g(u_n)\xrightarrow[n\to\infty]{C([0,T];L^\frac2m(\R^N))}g(u),
\end{align}
for any $T>0.$ In particular, $u$ solves~\eqref{nls} in $\Dr^\p\big((0,\infty)\times\R^N\big).$
\end{enumerate}
\medskip
\end{rmk}

\begin{thm}[\textbf{Existence and uniqueness of $\bs{L^2}$-solutions}]
\label{thmweak}
Let Assumption~$\ref{ass}$ be fulfilled and let $f\in L^1_\loc\big([0,\infty);L^2(\R^N)\big).$ Then for any $u_0\in L^2(\R^N),$ there exists a unique weak solution $u$ to \eqref{nls}--\eqref{u0}. In addition,
\begin{gather}
\label{Lm}
u\in L^{m+1}_\loc\big([0,\infty);L^{m+1}(\R^N)\big),	\\
\label{L2+}
\dfrac12\|u(t)\|_{L^2(\R^N)}^2+\Im(a)\dsp\vint_s^t\|u(\sigma)\|_{L^{m+1}(\R^N)}^{m+1}\d\sigma
			\le\dfrac12\|u(s)\|_{L^2(\R^N)}^2+\,\Im\dsp\iint\limits_{\;\;s\;\R^N}^{\text{}\;\;t}f(\sigma,x)\,\ovl{u(\sigma,x)}\,\d x\,\d\sigma,
\end{gather}
for any $t\ge s\ge0.$ Finally, if $v$ is a weak solution of \eqref{nls} with $v(0)=v_0\in L^2(\R^N)$ and $g\in L^1_\loc([0,\infty);L^2(\R^N))$ instead of $f$ in \eqref{nls} then,
\begin{gather}
\label{estthmweak}
\|u(t)-v(t)\|_{L^2(\R^N)}\le\|u(s)-v(s)\|_{L^2(\R^N)}+\vint_s^t\|f(\sigma)-g(\sigma)\|_{L^2(\R^N)}\d\sigma,
\end{gather}
for any $t\ge s\ge0.$
\end{thm}

\begin{rmk}
\label{rmkthmweak}
Let Assumption~$\ref{ass}$ be fulfilled. It follows from~\eqref{L2+} and Hölder's and Young's inequalities that if $f\in L^1\big((0,\infty);L^2(\R^N)\big)$ then,
\begin{gather*}
u\in L^\infty\big((0,\infty);L^2(\R^N)\big)\cap L^{m+1}\big((0,\infty);L^{m+1}(\R^N)\big).
\end{gather*}
By interpolation, we infer that for any $p\in[m+1,2),$
\begin{gather}
\label{rmkthmweak1}
u\in C_\b\big([0,\infty);L^2(\R^N)\big)\cap L^\frac{p(1-m)}{2-p}\big((0,\infty);L^p(\R^N)\big).
\end{gather}
If, in addition, $(\vphi_n)_{n\in\N}\subset L^2(\R^N),$ $(f_n)_{n\in\N}\subset L^1\big((0,\infty);L^2(\R^N)\big)$ and,
\begin{gather*}
\vphi_n\xrightarrow[n\to\infty]{L^2(\R^N)}u_0 \; \text{ and } \; f_n\xrightarrow[n\to\infty]{L^1((0,\infty);L^2(\R^N))}f,
\end{gather*}
then by~\eqref{estthmweak}, \eqref{rmkthmweak1} and again by interpolation, we have for any $p\in(m+1,2),$
\begin{gather*}
u_n\xrightarrow[n\to\infty]{C_\b([0,\infty);L^2(\R^N))\cap L^\frac{p(1-m)}{2-p}((0,\infty);L^p(\R^N))}u,
\end{gather*}
where for each $n\in\N,$ $u_n$ is the weak solution of \eqref{nls} with $u_n(0)=\vphi_n$ and $f_n$ instead of $f.$
\end{rmk}

\begin{thm}[\textbf{Existence and uniqueness of $\bs{H^1}$-solutions}]
\label{thmstrongH1}
Let Assumption~$\ref{ass}$ be fulfilled and let $f\in W^{1,1}_\loc\big([0,\infty);H^1(\R^N)\big).$ Then for any $u_0\in H^1(\R^N),$ there exists a unique $H^1$-solution $u$ to \eqref{nls}--\eqref{u0}. Furthermore, $u$ is also a weak solution and satisfies the following properties.
\begin{enumerate}[$1)$]
\item
\label{thmstrongH11}
$u\in C\big([0,\infty);L^2(\R^N)\big)\cap C^1\big([0,\infty);Y^\star\big)$ and $u$ satisfies \eqref{nls} in $Y^\star,$ for any $t\ge0,$ where $Y^\star=H^{-2}(\R^N)+L^\frac2m(\R^N).$
\item
\label{thmstrongH12}
$u\in C_\w\big([0,\infty);H^1(\R^N)\big)\cap W^{1,\infty}_\loc\big([0,\infty);H^{-1}(\R^N)+L^\frac2m(\R^N)\big)$ and,
\begin{gather}
\label{strongH11}
\|\nabla u(t)\|_{L^2(\R^N)}\le\|\nabla u_0\|_{L^2(\R^N)}+\int_0^t\|\nabla f(s)\|_{L^2(\R^N)}\d s,
\end{gather}
for any $t\ge0.$
\item
\label{thmstrongH13}
The map $t\longmapsto\|u(t)\|_{L^2(\R^N)}^2$ belongs to $W^{1,1}_\loc\big([0,\infty);\R\big)$ and we have,
\begin{gather}
\label{L2}
\frac12\frac{\d}{\d t}\|u(t)\|_{L^2(\R^N)}^2+\Im(a)\|u(t)\|_{L^{m+1}(\R^N)}^{m+1}=\Im\vint_{\R^N}f(t,x)\,\ovl{u(t,x)}\,\d x,
\end{gather}
for almost every $t>0.$
\end{enumerate}
\end{thm}

\begin{thm}[\textbf{Existence and uniqueness of $\bs{H^2}$-solutions}]
\label{thmstrongH2}
Let Assumption~$\ref{ass}$ be fulfilled and let $f\in W^{1,1}_\loc\big([0,\infty);L^2(\R^N)\big).$ Then for any $u_0\in H^2(\R^N)\cap L^{2m}(\R^N),$  there exists a unique $H^2$-solution $u$ to \eqref{nls}--\eqref{u0}. Furthermore, $u$ satisfies~\eqref{nls} in $L^2(\R^N),$ for almost every $t>0,$ and the following properties.
\begin{enumerate}[$1)$]
\item
\label{thmstrongH21}
$u\in C\big([0,\infty);H^1(\R^N)\cap L^{m+1}(\R^N)\big)\cap C^1\big([0,\infty);H^{-1}(\R^N)+L^\frac{m+1}m(\R^N)\big)$ and $u$ satisfies \eqref{nls} in $H^{-1}(\R^N)+L^\frac{m+1}m(\R^N),$ for any $t\ge0.$
\item
\label{thmstrongH22}
$u\in W^{1,\infty}_\loc\big([0,\infty);L^2(\R^N)\big)\cap L^\infty_\loc\big([0,\infty);H^2(\R^N)\cap L^{2m}(\R^N)\big)$ and,
\begin{empheq}[left=\empheqlbrace]{align}
\label{strongH21}
	&	\|u(t)-u(s)\|_{L^2(\R^N)}\le\|u_t\|_{L^\infty((s,t);L^2(\R^N))}|t-s|,			\dfrac{}{}				\\
\label{strongH22}
	&	\|\nabla u(t)-\nabla u(s)\|_{L^2(\R^N)}\le M|t-s|^\frac12,									\\
\label{strongH23}
	&	\left\|u_t\right\|_{L^\infty((0,t);L^2(\R^N))}
		\le\|\Delta u_0+a|u_0|^{m-1}u_0-f(0)\|_{L^2(\R^N)}+\int_0^t\|f^\p(\sigma)\|_{L^2(\R^N)}\d\sigma,
\end{empheq}
for any $t\ge s\ge0,$ where $M^2=2\|u_t\|_{L^\infty((s,t);L^2(\R^N))}\|\Delta u\|_{L^\infty((s,t);L^2(\R^N))}.$
\item
\label{thmstrongH23}
The map $t\longmapsto\|u(t)\|_{L^2(\R^N)}^2$ belongs to $C^1\big([0,\infty);\R\big)$ and~\eqref{L2} holds for any $t\ge0.$
\item
\label{thmstrongH24}
If $f\in W^{1,1}\big((0,\infty);L^2(\R^N)\big)$ then we have,
\begin{align*}
	&	u\in C_\b\big([0,\infty);H^1(\R^N)\big)\cap L^\infty\big((0,\infty);H^2(\R^N)\cap L^{2m}(\R^N)\big)\cap
			W^{1,\infty}\big((0,\infty);L^2(\R^N)\big).
\end{align*}
\end{enumerate}
\end{thm}

\begin{rmk}
\label{rmkf0}
Since $f\in W^{1,1}_\loc\big([0,\infty);L^2(\R^N)\big)\inj C\big([0,\infty);L^2(\R^N)\big)$ (see, for instance, 1) of Lemma~A.4 in Bégout and D\'iaz~\cite{MR4053613}), estimate~\eqref{strongH23} with $f(0)$ makes sense.
\end{rmk}

\begin{rmk}
\label{rmkthmstrong}
We recall that if $u\in L^2(\R^N)$ with $\Delta u\in L^2(\R^N)$ then $u\in H^2(\R^N).$ Furthermore, if $\|u\|^2_{H^{2,2}(\R^N)}=\|u\|_{L^2(\R^N)}^2+\|\Delta u\|_{L^2(\R^N)}^2$ then $\|\:.\:\|_{H^{2,2}(\R^N)}$ and $\|\:.\:\|_{H^2(\R^N)}$ are equivalent norms. Indeed, this us due to the Fourier transform and Plancherel's formula. Finally, note that,
\begin{gather}
\label{E}
\|\nabla u\|_{L^2(\R^N)}^2\le\|u\|_{L^2(\R^N)}\|\Delta u\|_{L^2(\R^N)}
\le\|u\|_{L^2(\R^N)}^2+\|\Delta u\|_{L^2(\R^N)}^2,
\end{gather}
for any $u\in H^2(\R^N).$
\end{rmk}

\begin{rmk}
\label{rmkthmexiuni}
Using a radically different method than the one we propose here, we may show that all the results of this section remain valid if we replace $\R^N$ with an unbounded domain $\Omega\neq\R^N.$ This will be the subject of a future work.
\end{rmk}

\section{Finite time extinction and asymptotic behavior}
\label{finite}

Following the method by Carles and Gallo~\cite{MR2765425} (also used by Carles and Ozawa~\cite{MR3306821}) and Bégout and D\'iaz~\cite{MR4053613}, we are able to prove the finite time extinction and asymptotic behavior results.

\begin{thm}
\label{thmextH2}
Let Assumption~$\ref{ass}$ be fulfilled with $N\in\{1,2,3\},$ let $f\in W^{1,1}\big((0,\infty);L^2(\R^N)\big),$ let $u_0\in H^1(\R^N)$ and assume that one of the following hypotheses holds.
\begin{enumerate}[$1)$]
\item
$N=1$ and $f\in W^{1,1}\big((0,\infty);H^1(\R)\big).$
\item
$N\in\{1,2,3\}$ and $u_0\in H^2(\R^N)\cap L^{2m}(\R^N).$ 
\end{enumerate}
Let $u$ be the unique strong solution of \eqref{nls}--\eqref{u0}. Finally, assume that there exists $T_0\ge0$ such that,
\begin{gather*}
\text{for almost every } \; t>T_0, \; f(t)=0.
\end{gather*}
Let $\ell$ be the exponant in $u_0\in H^\ell(\R^N).$  We have the following results.
\begin{enumerate}[$a)$]
\item
\label{thmextH21}
There exists a finite time $T_\star\ge T_0$ such that,
\begin{gather}
\label{0}
\forall t\ge T_\star, \; \|u(t)\|_{L^2(\R^N)}=0.
\end{gather}
Furthermore,
\begin{gather}
\label{T*}
T_\star\le C\|u\|_{L^\infty((0,\infty);H^\ell(\R^N))}^\frac{N(1-m)}{2\ell}\|u(T_0)\|_{L^2(\R^N)}^\frac{(1-m)(2\ell-N)}{2\ell}+T_0,
\end{gather}
where $C=C(\Im(a),N,m,\ell).$
\item
\label{thmextH22}
There exists $\eps_\star=\eps_\star(|a|,N,m)$ satisfying the following property. Let $\delta=\frac{(2\ell+N)+m(2\ell-N)}{4\ell}\in\left(\frac12,1\right).$ If $f\in W^{1,1}\big((0,\infty);H^1(\R^N)\big),$
\begin{gather*}
\begin{cases}
\left(\|u_0\|_{H^1(\R^N)}+\|f\|_{L^1((0,\infty);H^1(\R^N))}\right)^{1-m}\le\eps_\star\min\big\{1,T_0\big\},
&	\text{if } N=1,	\\
\left(\|u_0\|_{H^2(\R^N)}^m+\|f\|_{W^{1,1}((0,\infty);H^1(\R^N))}^m\right)^{1-m}\le\eps_\star\min\big\{1,T_0\big\},
&	\text{if } N\in\{2,3\},
\end{cases}
\end{gather*}
and if for almost every $t>0,$
\begin{gather}
\label{thmextH23}
\|f(t)\|_{L^2(\R^N)}^2\le\eps_\star\big(T_0-t\big)_+^\frac{2\delta-1}{1-\delta},
\end{gather}
then \eqref{0} holds with $T_\star=T_0.$
\end{enumerate}
\end{thm}

\begin{rmk}
\label{delta}
If $(N,\ell)\in\{(1,1),(2,2)\}$ then $\frac{2\delta-1}{1-\delta}=2\frac{1+m}{1-m},$ if $(N,\ell)=(1,2)$ then $\frac{2\delta-1}{1-\delta}=2\frac{1+3m}{3(1-m)}$ and if $(N,\ell)=(3,2)$ then $\frac{2\delta-1}{1-\delta}=2\frac{3+m}{1-m}.$ Note that if $N=1$ and $u_0\in H^2(\R^N)$ then there are two possible choices for $\frac{2\delta-1}{1-\delta}$ in~\eqref{thmextH23}: $2\frac{1+m}{1-m}$ or $2\frac{1+3m}{3(1-m)}.$ Since for $t$ near $T_0,$ $T_0-t<1$ then the choice the less restrictive is that for which $\frac{2\delta-1}{1-\delta}$ is the smallest as possible, that is $2\frac{1+3m}{3(1-m)}.$
\end{rmk}

\begin{rmk}
\label{rhmthmext}
In the case of our nonlinearity, Theorem~\ref{thmextH2} is an improvement of the result of Carles and Ozawa~\cite{MR3306821} in the sense they obtain the same conclusion as in $\ref{thmextH21})$ but with a presence harmonic confinement in~\eqref{nls}, $\Re(a)=0,$ $f=0,$ $N\in\{1,2\}$ and $\big(u_0\in H^1(\R)\cap\Fr(H^1(\R))\footnote{\label{fn3}$\Fr(H^1(\R))\inj L^{2m}(\R)$ and $\Fr(H^2(\R^2))\inj L^{2m}(\R^2),$ for any $\frac13<m\le1.$}\big),$ if $N=1$ and \big($u_0\in H^2(\R^2)\cap\Fr(H^2(\R^2))^{\ref{fn3}},$ $\|u_0\|_{L^2(\R^2)}$ small enough and $\frac12\le m<1\big),$ if $N=2.$ Additional nonlinearities are also considered in~\cite{MR3306821}.
\end{rmk}

\begin{thm}
\label{thm0s}
Let Assumption~$\ref{ass}$ be fulfilled with $N\ge4,$ let $f\in W^{1,1}_\loc\big([0,\infty);L^2(\R^N)\big)$ and let $u_0\in H^1(\R^N).$ Suppose further that $f\in W^{1,1}_\loc\big([0,\infty);H^1(\R^N)\big)$ or $u_0\in H^2(\R^N).$ Let $u$ be the unique strong solution of \eqref{nls}--\eqref{u0}. Finally, assume that there exists $T_0\ge0$ such that,
\begin{gather*}
\text{for almost every } t>T_0, \; f(t)=0.
\end{gather*}
Then we have for any $t\ge T_0,$
\begin{gather*}
\|u(t)\|_{L^2(\R^N)}\le\|u(T_0)\|_{L^2(\R^N)}e^{-C(t-T_0)},
\end{gather*}
if $N=4$ and $u_0\in H^2(\R^N),$
\begin{gather*}
\|u(t)\|_{L^2(\R^N)}\le\dfrac{\|u(T_0)\|_{L^2(\R^N)}}
{\left(1+C\|u(T_0)\|_{L^2(\R^N)}^\frac{(1-m)(N-2\ell)}{2\ell}(t-T_0)\right)^\frac{2\ell}{(1-m)(N-2\ell)}},
\end{gather*}
if $N\ge5$ or $u_0\in H^1(\R^N),$ where $C=C(\|u\|_{L^\infty((0,\infty);H^\ell(\R^N))},\Im(a),N,m,\ell).$
\end{thm}

\begin{thm}
\label{thm0w}
Let Assumption~$\ref{ass}$ be fulfilled, let $f\in L^1_\loc\big([0,\infty);L^2(\R^N)\big),$ let $u_0\in L^2(\R^N)$ and let $u$ be the unique weak solution of \eqref{nls}--\eqref{u0}. If
\begin{gather*}
f\in L^1\big((0,\infty);L^2(\R^N)\big),
\end{gather*}
then,
\begin{gather*}
\lim_{t\nearrow\infty}\|u(t)\|_{L^2(\R^N)}=0.
\end{gather*}
\end{thm}

\section{Proofs of the existence and uniqueness theorems}
\label{proofexi}

Since we have to prove existence in the whole space, the method is radically different than that used in Bégout and D\'iaz~\cite{MR4053613}.

\begin{thm}
\label{thmexicom}
Let Assumption~$\ref{ass}$ be fulfilled and let $\lambda,b_0>0.$ Then for any $F\in L^2(\R^N),$ there exists a unique solution $u$ to,
\begin{gather}
\label{nlss}
\begin{cases}
u\in H^2(\R^N)\cap L^{2m}(\R^N),						\medskip \\
-\lambda\Delta u-a\lambda|u|^{-(1-m)}u-\vi b_0u=F, \; \text{ in } \; L^2(\R^N).
\end{cases}
\end{gather}
In addition,
\begin{gather}
\label{est}
\|u\|_{H^2(\R^N)}^2+\|u\|_{L^{m+1}(\R^N)}^{m+1}+\|u\|_{L^{2m}(\R^N)}^{2m}\le M\|F\|_{L^2(\R^N)}^2,
\end{gather}
where $M=M(|a|,\Arg(a),b_0,\lambda).$ Furthermore, if $F$ is compactly supported then so is $u.$ Finally, let $G\in L^2(\R^N).$ If $v$ is a solution to~\eqref{nlss} with $G$ instead of $F$ then,
\begin{gather}
\label{uni}
\|u-v\|_{L^2(\R^N)}\le\frac1{b_0}\|F-G\|_{L^2(\R^N)}.
\end{gather}
Here and after, $\Arg(a)\in(0,\pi)$ denotes the principal value of the argument of $a.$
\end{thm}

\noindent
The proof of the theorem relies on the following lemmas.

\begin{lem}
\label{lemb}
Let Assumption~$\ref{ass}$ be fulfilled. Then there exists $b\in\C,$ with $|b|=1,$ satisfying the following property.
\begin{gather}
\label{b}
\Re(b)>0 \; \text{ and } \; \Im(b)<0,	\\
\label{ab}
2\sqrt m\,\Im(ab)>(1-m)\Re(ab)\ge0.
\end{gather}
In addition, $b=b(\Arg(a)).$ In particular, $ab$ satisfies \eqref{a}--\eqref{as} of Assumption~$\ref{ass}.$
\end{lem}

\begin{proof*}
Let $\theta_a=\Arg(a)\in(0,\pi),$ since $\Im(a)>0.$ We look for $b=e^{-\vi\theta_b},$ where $0<\theta_b<\frac\pi2.$ \\
\textbf{Case 1:} $\Re(a)<0.$ \\
If follows that, $\frac\pi2<\theta_a<\pi.$ We choose $\theta_b=\theta_a-\frac\pi2.$ We then have $ab=\vi|a|$ and the conclusion is clear. \\
\textbf{Case 2:} $\Re(a)\ge0.$ \\
If follows that, $0<\theta_a\le\frac\pi2$ and by \eqref{as}, one has
\begin{gather}
\label{demlemb}
2\sqrt m\sin(\theta_a)>(1-m)\cos(\theta_a)\ge0.
\end{gather}
By continuity and \eqref{demlemb}, there exists $\theta_b\in(0,\theta_a)$ such that,
\begin{gather}
2\sqrt m\sin(\theta_a-\theta_b)>(1-m)\cos(\theta_a-\theta_b)>0.
\end{gather}
Then, $0<\theta_a-\theta_b<\frac\pi2,$ $ab=|a|e^{\vi(\theta_a-\theta_b)}$ and again the conclusion is clear.
\medskip
\end{proof*}

\noindent
We may summarize the proof of Lemma~\ref{lemb} with the picture below.
\medskip

\begin{tabular}{rl}
\label{picture}
\begin{tikzpicture}[scale=1.1]
\fill[fill=gray!10] (0,0) -- (2.2,3.74) -- (-2.2,3.74);
\draw [->] (-2.5,0) -- (2.5,0);										
\draw [->] (0,-0.5) -- (0,4);											
\draw [domain=-2.2:0] plot (\x, {-1.7*\x});								
\draw [purple, domain=0:2.2] plot (\x, {1.7*\x});							
\draw (1,-0.1) -- (1,0.1);											
\draw (-0.1,1) -- (0.1,1);											
\draw (110:3) node[]{.};											
\draw (-20:1) node[]{.};											
\draw (-0.1,3) -- (0.1,3);											
\draw (0,0) -- (110:3);											
\draw (0,0) -- (-20:1);												
\draw [->] (110:2.5) arc (110:90:2.5);									
\draw [->] (0:0.8) arc (0:-20:0.8);									
\draw [->] (20:2.5) arc (20:50:2.5);									
\draw (0,0) node[below left]{$0$};									
\draw (1,0) node[above right]{$1$};									
\draw (0,1) node[above right]{$\vi$};									
\draw (109:3) node[above]{$a=|a|e^{\vi\theta_a}$};						
\draw (-20:1) node[below right]{$b=e^{-\vi\theta_b}$};					
\draw (0,3) node[above right]{$ab$};									
\draw (35:2.7) node[]{$+$};										
\draw (100:2.8) node[]{$-\theta_b$};									
\draw (-15:0.75) node[right]{$\longleftarrow-\theta_b$};					
\draw (2.55,0) node[right]{$\Re(z)$};									
\draw (0,4) node[above]{$\Im(z)$};									
\draw (1.55,2.5) node[right]{$\scriptstyle\Im(z)=\frac{1-m}{2\sqrt m}|\Re(z)|$};	
\draw (0,-1.3) node[]{$\theta_b=\theta_a-\frac\pi2$};
\draw (0,-1.7) node[]{Case 1: $\Re(a)<0$};
\end{tikzpicture}
&
\begin{tikzpicture}[scale=1.1]
\fill[fill=gray!10] (0,0) -- (2.2,3.74) -- (-2.2,3.74);
\draw [->] (-2.5,0) -- (2.5,0);										
\draw [->] (0,-0.5) -- (0,4);											
\draw [domain=-2.2:0] plot (\x, {-1.7*\x});								
\draw [purple, domain=0:2.2] plot (\x, {1.7*\x});							
\draw (1,-0.1) -- (1,0.1);											
\draw (-0.1,1) -- (0.1,1);											
\draw (85:3) node[]{.};											
\draw (-20:1) node[]{.};											
\draw (65:3) node[]{.};											
\draw (0,0) -- (85:3);												
\draw (0,0) -- (-20:1);												
\draw (0,0) -- (65:3);												
\draw [->] (85:2) arc (85:65:2);										
\draw [->] (0:0.8) arc (0:-20:0.8);									
\draw [->] (20:2.5) arc (20:50:2.5);									
\draw (0,0) node[below left]{$0$};									
\draw (1,0) node[above right]{$1$};									
\draw (0,1) node[above left]{$\vi$};									
\draw (85:3) node[above]{$a$};										
\draw (-20:1) node[below right]{$b=e^{-\vi\theta_b}$};					
\draw (65:3) node[above]{$ab$};									
\draw (35:2.7) node[]{$+$};										
\draw (75:2.3) node[]{$-\theta_b$};									
\draw (-15:0.75) node[right]{$\longleftarrow-\theta_b$};					
\draw (2.55,0) node[right]{$\Re(z)$};									
\draw (0,4) node[above]{$\Im(z)$};									
\draw (1.55,2.5) node[right]{$\scriptstyle\Im(z)=\frac{1-m}{2\sqrt m}|\Re(z)|$};	
\draw (0,-1.3) node[]{$0<\theta_b\ll1$};
\draw (0,-1.7) node[]{Case 2: $\Re(a)\ge0$};
\end{tikzpicture}
\medskip
\end{tabular}

\begin{lem}
\label{lemmon}
Let $0<m<1.$ Set for any $z\in\C,$ $g(z)=|z|^{-(1-m)}z$ $(g(0)=0).$ We define the mapping for any measurable function $u:\R^N\tends\C,$ which we still denote by $g,$ by $g(u)(x)=g(u(x)).$ Then for any $p\in[1,\infty),$
\begin{gather}
\label{lemmon1}
g\in C\big(L^p(\R^N);L^\frac{p}m(\R^N)\big) \text{ and } g \text{ is bounded on bounded sets.}
\end{gather}
Let $a\in\C$ with $\Im(a)>0$ satisfying~\eqref{a}. Then $\big(g(u)-g(v)\big)(\ovl{u-v})\in L^1(\R^N)$ and,
\begin{gather}
\label{lemmon2}
\Re\left(-\vi\,a\vint_{\R^N}\big(g(u)-g(v)\big)(\ovl{u-v})\d x\right)\ge0,
\end{gather}
for any $u,v\in L^{m+1}(\R^N).$
\end{lem}

\begin{proof*}
Property~\eqref{lemmon1} is an obvious consequence of \eqref{g} which implies the integrability property in the lemma. By Lemma~2.2 of Liskevich and Perel$^\p$muter~\cite{MR1224619}, we have
\begin{gather}
\label{demlemmon}
2\sqrt m\left|\Im\Big(\big(g(z_1)-g(z_2)\big)\big(\ovl{z_1-z_2}\big)\Big)\right|
		\le(1-m)\Re\Big(\big(g(z_1)-g(z_2)\big)\big(\ovl{z_1-z_2}\big)\Big),
\end{gather}
for any $(z_1,z_2)\in\C^2.$ Let $u,v\in L^{m+1}(\R^N).$ We have by \eqref{demlemmon},
\begin{align*}
	& \; \Re\left(-\vi\,a\vint_{\R^N}\big(g(u)-g(v)\big)(\ovl{u-v})\d x\right)								\\
   =	& \; \Im(a)\Re\vint_{\R^N}\big(g(u)-g(v)\big)\big(\ovl{u-v}\big)\d x
   		+\Re(a)\Im\vint_{\R^N}\big(g(u)-g(v)\big)\big(\ovl{u-v}\big)\d x								\\
 \ge	& \; \left(\Im(a)-|\Re(a)|\frac{1-m}{2\sqrt m}\right)\Re\vint_{\R^N}\big(g(u)-g(v)\big)\big(\ovl{u-v}\big)\d x	\\
 \ge	& \; 0.
\end{align*}
The lemma is proved.
\medskip
\end{proof*}

\begin{lem}[\textbf{\cite{MR4053613}}]
\label{lemLap}
Let $0<m<1$ and let $a\in\C$ with $\Im(a)>0$ satisfying~\eqref{a}. Let $g$ be as in Lemma~$\ref{lemmon}.$ Then $g(u)\ovl{\Delta u}\in L^1(\R^N)$ and,
\begin{gather}
\label{lemLap1}
\Re\left(\vi a\vint_{\R^N}g(u)\ovl{\Delta u}\d x\right)\ge0,
\end{gather}
for any $u,v\in H^2(\R^N)\cap L^{2m}(\R^N).$
\end{lem}

\begin{proof*}
See Bégout and D\'iaz~\cite{MR4053613} (Lemma~6.3).
\medskip
\end{proof*}

\begin{vproof}{of Theorem~\ref{thmexicom}.}
Let Assumption~\ref{ass} be fulfilled, $\lambda,b_0>0$ and $F\in L^2(\R^N).$ Let $g$ be as in Lemma~\ref{lemmon}. We want to solve,
\begin{gather}
\label{uF}
-\lambda\Delta u-a\lambda g(u)-\vi b_0u=F, \; \text{ in } \; H^{-1}(\R^N)+L^\frac{m+1}m(\R^N).	\tag{$u_F$}
\end{gather}
We proceed with the proof in five steps. \\
\textbf{Step~1: A first estimate.} Let $G\in L^2(\R^N).$ If $u,v\in H^2_\loc(\R^N)\cap H^1(\R^N)\cap L^{m+1}(\R^N)$ are solutions of $(u_F)$ and $(v_G),$ respectively, then estimate \eqref{uni} holds true. \\
We multiply by $\vi\ovl{\vphi},$ for $\vphi\in\Dr(\R^N),$ the equation satisfied by $u-v,$ we integrate by parts and we take the real part. By density of $\Dr(\R^N)$ in $H^1(\R^N)\cap L^{m+1}(\R^N)$ and \eqref{lemmon1}, $\big(g(u)-g(v)\big)(\ovl{u-v})\in L^1(\R^N)$ and we may choose $\vphi=u-v.$ It follows that,
\begin{gather}
\label{demthmexicom1}
\lambda\Re\left(-\vi a\vint_{\R^N}\big(g(u)-g(v)\big)(\ovl{u-v})\d x\right)+b_0\|u-v\|_{L^2(\R^N)}^2
=-\Im\left(\;\vint_{\R^N}(F-G)\ovl{(u-v)}\d x\right).
\end{gather}
Estimate \eqref{uni} then comes from \eqref{demthmexicom1}, \eqref{lemmon2} and Cauchy-Schwarz's inequality. \\
\textbf{Step~2: A second estimate.} If $u$ is a solution to \eqref{nlss} then $u\in L^{m+1}(\R^N)$ and satisfies \eqref{est}. \\
Since $2m<m+1<2,$ then $L^{2m}(\R^N)\cap L^2(\R^N)\subset L^{m+1}(\R^N).$ By Theorem~2.9 in Bégout and D\'iaz~\cite{MR3315701},
\begin{gather}
\label{demthmexicom3}
\|u\|_{H^1(\R^N)}^2+\|u\|_{L^{m+1}(\R^N)}^{m+1}\le M(|a|,b_0,\lambda)\|F\|_{L^2(\R^N)}^2.
\end{gather}
Let $b\in\C$ be given by Lemma~\ref{lemb}. We multiply the equation in \eqref{nlss} by $-\vi b\ovl{\Delta u},$ integrate by parts and take the real part. We obtain,
\begin{equation}
\begin{aligned}
\label{demthmexicom5}
	&-\lambda\Im(b)\|\Delta u\|_{L^2(\R^N)}^2+\lambda\Re\left(\vi ab\vint_{\R^N}g(u)\ovl{\Delta u}\d x\right)+b_0\Re(b)\|\nabla u\|_{L^2(\R^N)}^2	\\
  =	&	\; \Im\left(b\int_{\R^N}F\ovl{\Delta u}\d x\right).
\end{aligned}
\end{equation}
By \eqref{ab}, we may apply Lemma~\ref{lemLap}. Using \eqref{b}, \eqref{lemLap1} and applying Cauchy-Schwarz's inequality in \eqref{demthmexicom5}, one obtains,
\begin{gather}
\label{demthmexicom6}
\|\Delta u\|_{L^2(\R^N)}\le\frac{|b|}{\lambda|\Im(b)|}\|F\|_{L^2(\R^N)}.
\end{gather}
Now, since by Plancherel's formula, $\|u\|_{\dot{H^2}(\R^N)}\le C\||\xi|^2\widehat u\|_{L^2(\R^N)}\le C\|\Delta u\|_{L^2(\R^N)},$ putting together \eqref{demthmexicom3} and \eqref{demthmexicom6}, one obtains \eqref{est}. \\
\textbf{Step~3: Compactness of the solution.} If $\supp F$ is compact and if $u\in H^1(\R^N)\cap L^{m+1}(\R^N)$ is a solution to \eqref{uF} then $\supp u$ is compact. \\
This comes from Theorem~3.6 in Bégout and D\'iaz~\cite{MR2876246}. \\
\textbf{Step~4: Existence and uniqueness.} There exists a unique solution $u\in H^2_\loc(\R^N)\cap H^1(\R^N)\cap L^{m+1}(\R^N)$ to \eqref{uF}. \\
By Theorem~2.8 in Bégout and D\'iaz~\cite{MR3315701}, equation \eqref{uF} admits a solution $u\in H^1(\R^N)\cap L^{m+1}(\R^N).$ By Proposition~4.5 in Bégout and D\'iaz~\cite{MR2876246}, $u\in H^2_\loc(\R^N).$ Finally, by Step~1 this solution is unique. \\
\textbf{Step~5: Conclusion.} \\
Estimates \eqref{est}--\eqref{uni}, uniqueness and compactness property come from Steps~1--3, once the existence of a solution to \eqref{nlss} is proved. Let $u\in H^2_\loc(\R^N)\cap H^1(\R^N)\cap L^{m+1}(\R^N)$ the solution of \eqref{uF} be given by Step~4. Let $(F_n)_{n\in\N}\subset\Dr(\R^N)$ be such that $F_n\xrightarrow[n\to\infty]{L^2(\R^N)}F.$ Finally, for each $n\in\N,$ denote by $u_n$ the unique solution to \eqref{nlss}, where the right-hand side is $F_n$ instead of $F$ (Steps~4 and 3). By Steps~1 and 2, $(u_n)_{n\in\N}$ is bounded in $H^2(\R^N)$ and $u_n\xrightarrow[n\to\infty]{L^2(\R^N)}u.$ It follows that $u\in H^2(\R^N)$ and, from the equation in \eqref{nlss}, $g(u)\in L^2(\R^N).$ Hence $u$ is a solution to \eqref{nlss}. This concludes the proof of the lemma.
\medskip
\end{vproof}

\begin{cor}
\label{Amaxmon}
Let Assumption~$\ref{ass}$ be fulfilled. Let us define the following $($nonlinear$)$ operator on $L^2(\R^N).$
\begin{gather*}
\begin{cases}
D(A)=H^2(\R^N)\cap L^{2m}(\R^N),				\medskip \\
\forall u\in D(A), \;  Au=-\vi\Delta u-\vi a|u|^{-(1-m)}u,
\end{cases}
\end{gather*}
Then $A$ is maximal monotone on $L^2(\R^N)$ $($and so $m$-accretive$)$ with dense domain.
\end{cor}

\begin{proof*}
The density is obvious. For any $\lambda>0,$ $I+\lambda A$ is bijective from $D(A)$ onto $L^2(\R^N)$ and $(I+\lambda A)^{-1}$ is a contraction (Theorem~\ref{thmexicom}). It follows that $A$ is maximal monotone (Brezis~\cite{MR0348562}, Proposition~2.2, p.23).
\medskip
\end{proof*}

\begin{vproof}{of Theorem~\ref{thmstrongH2}.}
Let $g$ be as in Lemma~\ref{lemmon}. We first recall that by Remark~\ref{rmkf0},
\begin{gather}
\label{demH20}
f\in C\big([0,\infty);L^2(\R^N)\big).
\end{gather}
By Corollary~\ref{Amaxmon} and Barbu~\cite{MR0390843} (Theorem~2.2, p.131), there exists a unique $u\in W^{1,\infty}_\loc\big([0,\infty);L^2(\R^N)\big)$ satisfying $u(t)\in H^2(\R^N)\cap L^{2m}(\R^N)$ and \eqref{nls} in $L^2(\R^N),$ for almost every $t>0,$ $u(0)=u_0$ and \eqref{strongH23}. This last estimate yields \eqref{strongH21}. Since $u\in W^{1,\infty}_\loc\big([0,\infty);L^2(\R^N)\big),$ it follows from Lemma~A.5 in Bégout and D\'iaz~\cite{MR4053613} that the map $M:t\longmapsto\frac12\|u(t)\|_{L^2(\R^N)}^2$ belongs to $W^{1,\infty}_\loc\big([0,\infty);\R\big)$ and $M^\p(t)=\big(u(t),u_t(t)\big)_{L^2(\R^N)},$ for almost every $t>0.$ Multiplying \eqref{nls} by $\ovl{\vi u},$ integrating by parts over $\R^N$ and taking the real part, we obtain \eqref{L2}, for almost every $t>0.$ We deduce easily from \eqref{L2}, \eqref{demH20} and Hölder's inequality that $u\in L^\infty_\loc\big([0,\infty);L^{m+1}(\R^N)\big).$ Multiplying again \eqref{nls} by $\ovl u,$ integrating by parts and taking the real part, we get
\begin{gather*}
\|\nabla u(t)\|_{L^2(\R^N)}^2
\le|\Re(a)|\|u(t)\|_{L^{m+1}(\R^N)}^{m+1}+\left(\|u_t(t)\|_{L^2(\R^N)}+\|f(t)\|_{L^2(\R^N)}\right)\|u(t)\|_{L^2(\R^N)},
\end{gather*}
for almost every $t>0.$ It follows that $u\in L^\infty_\loc\big([0,\infty);H^1(\R^N)\big).$ We infer that $u$ is an $H^2$-solution. Let $b\in\C$ be given by Lemma~\ref{lemb}. We multiply \eqref{nls} by $\ovl{\vi abg(u)},$ integrate and take the real part. We get,
\begin{equation}
\begin{aligned}
\label{demH21}
	&	\Re\left(\ovl{ab}\vint_{\R^N}u_t\ovl{g(u)}\d x\right)+\Re\left(\ovl{\vi ab}\vint_{\R^N}\ovl{g(u)}\Delta u\d x\right)
			+|a|^2\Re(\ovl{\vi b})\|g(u)\|_{L^2(\R^N)}^2		\\
   =	&	\; \Re\left(\ovl{\vi ab}\vint_{\R^N}f\ovl{g(u)}\d x\right).
\end{aligned}
\end{equation}
By Lemma~\ref{lemb}, we have \eqref{lemLap1}. This implies,
\begin{gather}
\label{demH22}
\Re\left(\ovl{\vi ab}\vint_{\R^N}\ovl{g(u)}\Delta u\d x\right)=\Re\left(\vi ab\vint_{\R^N}g(u)\ovl{\Delta u}\d x\right)\ge0,
\end{gather}
and \eqref{demH21} becomes,
\begin{gather}
\label{demH23}
|a||\Im(b)|\,\|u\|_{L^{2m}(\R^N)}^{2m}\le\vint_{\R^N}|(u_t+\vi f)\ovl{g(u)}|\d x,
\end{gather}
since $\Re(\ovl{\vi b})=-\Im(b)>0,$ by \eqref{b}. By Cauchy-Schwarz's and Young's inequalities, we get
\begin{gather}
\label{demH24}
\vint_{\R^N}|(u_t+\vi f)\ovl{g(u)}|\d x\le\frac1{2|a||\Im(b)|}\|u_t+\vi f\|_{L^2(\R^N)}^2+\frac{|a||\Im(b)|}2\|u\|_{L^{2m}(\R^N)}^{2m}.
\end{gather}
Putting together \eqref{demH23} and \eqref{demH24}, we arrive at,
\begin{gather}
\label{demH25}
\|u(t)\|_{L^{2m}(\R^N)}^{2m}\le\frac1{|a|^2|\Im(b)|^2}\left(\|u_t(t)\|_{L^2(\R^N)}+\|f(t)\|_{L^2(\R^N)}\right)^2,
\end{gather}
for almost every $t>0.$ Multiplying again \eqref{nls} by $\vi b\ovl{\Delta u},$ using \eqref{demH22} and proceeding as above, we arrive at,
\begin{gather}
\label{demH26}
\|\Delta u(t)\|_{L^2(\R^N)}\le\frac1{|\Im(b)|}\left(\|u_t(t)\|_{L^2(\R^N)}+\|f(t)\|_{L^2(\R^N)}\right),
\end{gather}
for almost every $t>0.$
By \eqref{demH20}, \eqref{demH25}, \eqref{demH26}, Remark~\ref{rmkthmstrong} and Hölder's inequality (recalling that $2m<m+1<2),$ we obtain,
\begin{gather}
\label{demH27}
u\in L^\infty_\loc\big([0,\infty);H^2(\R^N)\big)\cap L^\infty_\loc\big([0,\infty);L^{2m}(\R^N)\big),	\\
\label{demH28}
u\in C\big([0,\infty);L^2(\R^N)\big)\cap  L^\infty_\loc\big([0,\infty);L^{2m}(\R^N)\big)\inj C\big([0,\infty);L^{m+1}(\R^N)\big).
\end{gather}
Recalling that $u\in W^{1,\infty}_\loc\big([0,\infty);L^2(\R^N)\big),$ by~\eqref{demH27} and the embedding 3) of Lemma~A.4, we have $u\in C\big([0,\infty);H^1(\R^N)\big).$ We then deduce Property~\ref{thmstrongH21}), with help of \eqref{gmm}, \eqref{demH20} and \eqref{nls}. With \eqref{E}, \eqref{strongH21} and \eqref{demH27}, we get \eqref{strongH22} and Property~\ref{thmstrongH22}) is proved. Property~\ref{thmstrongH23}) comes from \eqref{L2}, \eqref{demH20} and \eqref{demH28}. Finally, Property~\ref{thmstrongH24}) follows easily from Remarks~\ref{rmkthmweak}, \ref{rmkf0} and \ref{rmkthmstrong}, \eqref{strongH23}, \eqref{demH25} and \eqref{demH26}. This concludes the proof of the theorem.
\medskip
\end{vproof}

\begin{lem}
\label{lemdep}
Let Assumption~$\ref{ass}$ be fulfilled and $f,g\in L^1_\loc\big([0,\infty);L^2(\R^N)\big).$ If $u$ and $v$ are strong solutions or weak solutions of
\begin{gather*}
\vi u_t+\Delta u+a|u|^{-(1-m)}u=f_1,	\\
\vi v_t+\Delta v+a|v|^{-(1-m)}v=f_2,
\end{gather*}
respectively, then $u,v\in C\big([0,\infty);L^2(\Omega)\big)$ and
\begin{gather}
\label{lemdep1}
\|u(t)-v(t)\|_{L^2(\Omega)}\le\|u(s)-v(s)\|_{L^2(\Omega)}+\vint_s^t\|f_1(\sigma)-f_2(\sigma)\|_{L^2(\Omega)}\d\sigma,
\end{gather}
for any $t\ge s\ge0.$
\end{lem}

\begin{proof*}
Let $X=H^1(\R^N)\cap L^{m+1}(\R^N)$ and let $u,v$ be as in the lemma. Continuity comes from \eqref{Xcon} and Definition~\ref{defsol}. Estimate~\eqref{lemdep1} being stable by passing to the limit in $C\big([0,T];L^2(\R^N)\big)\times L^1\big((0,T);L^2(\R^N)\big),$ for any $T>0,$ it is sufficient to establish it for the $H^2$-solutions. And since an $H^2$-solution is an $H^1$ solution, we may assume that $u,v$ are $H^1$ solution. Making the difference between the two equations, it follows from \ref{rmkdefsol3}) of Remark~\ref{rmkdefsol} that we can take the $X^\star-X$ duality product of the result with $\vi(u-v).$ With help of (A.3) of Lemma~A.5 in Bégout and D\'iaz~\cite{MR4053613}, \eqref{dualg}, \eqref{lemmon2} and Cauchy-Schwarz's inequality, we then arrive at,
\begin{gather*}
\frac12\frac\d{\d t}\|u(\:.\:)-v(\:.\:)\|_{L^2(\Omega)}^2\le\|f_1-f_2\|_{L^2(\Omega)}\|u-v\|_{L^2(\Omega)},
\end{gather*}
almost everywhere on $(0,\infty).$ Integrating over $(s,t),$ one obtains \eqref{lemdep1}.
\medskip
\end{proof*}

\begin{vproof}{of Theorem~\ref{thmweak}.}
Existence, estimate~\eqref{estthmweak} and uniqueness comes from density of $H^2(\R^N)\times W^{1,1}_\loc([0,\infty);L^2(\R^N))$ in $L^2(\R^N)\times L^1_\loc([0,\infty);L^2(\R^N)),$ Theorem~\ref{thmstrongH2}, Lemma~\ref{lemdep} and completeness of $C\big([0,T];L^2(\R^N)\big),$ for any $T>0.$ Finally, estimates \eqref{Lm}--\eqref{L2+} are due to Bégout and D\'iaz~\cite{MR4053613} (Proposition~2.3). This ends the proof of the theorem.
\medskip
\end{vproof}

\begin{vproof}{of Theorem~\ref{thmstrongH1}.}
Uniqueness comes from Lemma~\ref{lemdep}. Let $f\in W^{1,1}_\loc([0,\infty);H^1(\R^N))$ and let $u_0\in H^1(\R^N).$ Let $(\vphi_n)_{n\in\N}\subset\Dr(\R^N)$ be such that $\vphi_n\xrightarrow[n\to\infty]{H^1(\R^N)}u_0.$ Finally, let $g$ be defined as in Lemma~\ref{lemmon} and for each $n\in\N,$ let $u_n$ the unique $H^2$-solution of \eqref{nls} such that $u_n(0)=\vphi_n,$ be given by Theorem~\ref{thmstrongH2}. By Lemma~\ref{lemdep}, we have for any $T>0$ and $n,p\in\N,$
\begin{gather}
\label{demthmstrongH11}
\|u_n\|_{C([0,T];L^2(\R^N))}\le\|\vphi_n\|_{L^2(\R^N)}+\int_0^T\|f(t)\|_{L^2(\R^N)}\d t,	\\
\nonumber
\|u_n-u_p\|_{L^\infty((0,\infty);L^2(\R^N))}\le\|\vphi_n-\vphi_p\|_{L^2(\R^N)},
\end{gather}
It follows that for any $T>0,$ $(u_n)_{n\in\N}$ is a Cauchy sequence in $C\big([0,T];L^2(\R^N)\big).$ As a consequence, there exists $u\in C\big([0,\infty);L^2(\R^N)\big)$ such that for any $T>0,$
\begin{gather}
\label{demthmstrongH12}
u_n\xrightarrow[n\to\infty]{C([0,T];L^2(\R^N))}u.
\end{gather}
By definition, it follows from \eqref{demthmstrongH12} that $u$ is a weak solution of \eqref{nls}--\eqref{u0}. By Theorem~\ref{thmstrongH2}, we can take the $L^2$-scalar product of \eqref{nls} with $-\vi\Delta u_n$ and it follows from (A.4) in Bégout and D\'iaz~\cite{MR4053613} that for any $n\in\N$ and almost every $s>0,$
\begin{gather*}
\frac12\frac{\d}{\d t}\|\nabla u_n(s)\|_{L^2(\R^N)}^2+\Re\left(\vi a\vint_{\R^N}g(u_n(s))\ovl{\Delta u_n(s)}\d x\right)
=\big(\nabla f(s),\vi\nabla u_n(s)\big)_{L^2(\R^N)}.
\end{gather*}
which gives with \eqref{lemLap1} and Cauchy-Schwarz's inequality,
\begin{gather*}
\frac12\frac{\d}{\d t}\|\nabla u_n(s)\|_{L^2(\R^N)}^2\le\|\nabla f(s)\|_{L^2(\R^N)}\|\nabla u_n(s)\|_{L^2(\R^N)}.
\end{gather*}
By integration, we obtain for any $t>0$ and any $n\in\N,$ 
\begin{gather}
\label{demthmstrongH15}
\|\nabla u_n(t)\|_{L^2(\R^N)}\le\|\nabla \vphi_n\|_{L^2(\R^N)}+\int_0^t\|\nabla f(s)\|_{L^2(\R^N)}\d s.
\end{gather}
By the Sobolev embedding (see, for instance, 1) of Lemma~A.4 in Bégout and D\'iaz~\cite{MR4053613}),
\begin{gather}
\label{demthmstrongH16}
W^{1,1}_\loc\big([0,\infty);L^2(\R^N)\big)\inj C\big([0,\infty);L^2(\R^N)\big),
\end{gather}
\eqref{demthmstrongH11}, \eqref{demthmstrongH15}, \eqref{lemmon1} and \eqref{nls}, we infer that,
\begin{gather}
\label{demthmstrongH17}
(u_n)_{n\in\N} \text{ is bounded in } L^\infty\big((0,T);H^1(\R^N)\big)\cap W^{1,\infty}\big((0,T);Z^\star\big),
\end{gather}
for any $T>0,$ where $Z^\star=H^{-1}(\R^N)+L^\frac2m(\R^N)$ is the topological dual space of $Z=H^1(\R^N)\cap L^\frac2{2-m}(\R^N).$ Note that $Z^\star$ is reflexive (Lemma~A.2 in Bégout and D\'iaz~\cite{MR4053613}) and since $H^1(\R^N)\inj Z^\star,$ it follows from \eqref{demthmstrongH12}, \eqref{demthmstrongH17}, \eqref{rmkdefsol41} and Proposition~1.1.2, p.2, and (ii) of Remark~1.3.13, p.12, in Cazenave~\cite{MR2002047} that,
\begin{align}
\label{demthmstrongH19}
&	u\in C_\w\left([0,\infty);H^1(\R^N)\right)\cap W^{1,\infty}_\loc\big([0,\infty);Z^\star\big),	\\
\label{demthmstrongH110}
&	\Delta u\in C\left([0,\infty);H^{-2}(\R^N)\right),									\\
\label{demthmstrongH111}
&	u_n(t)\weak u(t), \; \text{ in } \; H^1_\w(\R^N), \; \text{ as } \; n\to\infty,
\end{align}
for any $t\ge0.$ After integration of \eqref{L2}, we see with help of \eqref{demthmstrongH11} that for any $T>0,$ $(u_n)_{n\in\N}$ is bounded in $L^{m+1}\big((0,T);L^{m+1}(\R^N)\big)\cong L^{m+1}\big((0,T)\times\R^N\big),$ which is reflexive.
We infer with \eqref{demthmstrongH12},
\begin{gather}
\label{demthmstrongH113}
u\in L^{m+1}_\loc\big([0,\infty);L^{m+1}(\R^N)\big).
\end{gather}
By \ref{rmkdefsol4}) of Remark \ref{rmkdefsol}, \eqref{demthmstrongH16}, \eqref{demthmstrongH19}, \eqref{demthmstrongH113} and \eqref{nls}, it follows that $u$ satisfies \ref{defsol1}) of Definition~\ref{defsol} and then $u$ is an $H^1$-solution. By \ref{rmkdefsol3}) of Remark \ref{rmkdefsol}, we can take the $X-X^\star$ duality product with $\vi u,$ where $X=H^1(\R^N)\cap L^{m+1}(\R^N).$ Applying Lemma~A.5 of Bégout and D\'iaz~\cite{MR4053613} and \eqref{dualg}, Property~\ref{thmstrongH13}) follows. Estimate \eqref{strongH11} comes from \eqref{demthmstrongH111}, \eqref{demthmstrongH15} and the weak lower semicontinuity of the norm. Finally, smoothness of the solution in Properties \ref{thmstrongH11}) and \ref{thmstrongH12}) follows easily from \eqref{demthmstrongH16}, \eqref{demthmstrongH19}, \eqref{demthmstrongH110}, \eqref{lemmon1} and the equation \eqref{nls}. This concludes the proof of the theorem.
\medskip
\end{vproof}

\section{Proofs of the finite time extinction and asymptotic behavior theorems}
\label{proofext}

\begin{vproof}{of Theorem~\ref{thmextH2}.}
Apply Theorems~\ref{thmstrongH1}, \ref{thmstrongH2} and use the general theorem of finite time extinction in \cite{MR4053613} (Theorem~2.1 and Remark~4.8). Nevertheless, to make the proof more understandable, we briefly explain how to obtain \eqref{0}--\eqref{T*}. Let $\ell=1,$ if $u_0\in H^1(\R^N)$ and $\ell=2,$ if $u_0\in H^2(\R^N).$ Assume that for some $T_0\ge0,$ $f(t)=0,$ for almost every $t>T_0.$ It follows from Theorems~\ref{thmstrongH1}, \ref{thmstrongH2} and Remark~\ref{rmkthmweak} that $u\in L^\infty\big((0,\infty);H^\ell(\R^N)\big).$ We have by Gagliardo-Nirenberg's inequality and \eqref{L2},
\begin{gather*}
\|u(t)\|_{L^2(\R^N)}^\frac{(2\ell+N)+m(2\ell-N)}{2\ell}
\le C\|u\|_{L^\infty((0,\infty);H^\ell(\R^N))}^\frac{N(1-m)}{2\ell}\|u(t)\|_{L^{m+1}(\R^N)}^{m+1},	\\
\dfrac{\d}{\d t}\|u(t)\|_{L^2(\R^N)}^2+2\Im(a)\|u(t)\|_{L^{m+1}(\R^N)}^{m+1}=0,
\end{gather*}
for almost every $t>T_0.$ It follows that,
\begin{gather}
\label{demthmextH2}
y^\p(t)+C y(t)^\delta\le0,
\end{gather}
for almost every $t>T_0,$ where $y(t)=\|u(t)\|_{L^2(\R^N)}^2$ and $\delta=\frac{(2\ell+N)+m(2\ell-N)}{4\ell}.$ By our assumption on $\ell,$ we have $\delta\in(0,1)$ if $N\le3.$ Hence \eqref{0}--\eqref{T*} by integration.
\medskip
\end{vproof}

\begin{vproof}{of Theorem~\ref{thm0s}.}
Let $\ell=1,$ if $u_0\in H^1(\R^N)$ and $\ell=2,$ if $u_0\in H^2(\R^N).$ By Theorems~\ref{thmstrongH1}, \ref{thmstrongH2} and Remark~\ref{rmkthmweak}, $u\in L^\infty\big((0,\infty);H^\ell(\R^N)\big).$ Repeating the proof of Theorem \ref{thmextH2}, we obtain obtain \eqref{demthmextH2}. According to the different cases as in the theorem, we have $\delta=1$ or $\delta>1.$ The results then follow by integration (see also \eqref{intro2} and the lines below). For more details, see 3) of Remark~2.4 in \cite{MR4053613}.
\medskip
\end{vproof}

\begin{vproof}{of Theorem~\ref{thm0w}.}
By Remark~\ref{rmkthmweak}, we may assume that $f\in\Dr\big([0,\infty);L^2(\R^N)\big)$ and $u_0\in H^2(\R^N).$ Let $[0,T_0]\supset\supp f.$ By \eqref{L2}, $\frac{\d}{\d t}\|u(t)\|_{L^2(\R^N)}^2\le0,$ for any $t>T_0.$ It follows that
$\vlim_{t\nearrow\infty}\|u(t)\|_{L^2(\R^N)}=\ell_0,$ for some $\ell_0\in[0,\infty).$ Let $q\in(2,\infty)$ with $(N-2)q<2N.$ By Hölder's inequality and Sobolev's embedding $H^1(\R^N)\inj L^q(\R^N),$ there exists $\theta\in(0,1)$ such that,
\begin{gather*}
\ell_0\le\|u(t)\|_{L^2(\R^N)}\le\|u(t)\|_{L^{m+1}(\R^N)}^\theta\|u(t)\|_{L^q(\R^N)}^{1-\theta}\le C\|u(t)\|_{L^{m+1}(\R^N)}^\theta\|u\|_{L^\infty((0,\infty);H^1(\R^N))}^{1-\theta},
\end{gather*}
for any $t>T_0.$ We get, still by \eqref{L2},
\begin{gather*}
\frac{\d}{\d t}\|u(t)\|_{L^2(\R^N)}^2\le-C\ell_0^\frac{m+1}\theta\le0,
\end{gather*}
for any $t>T_0.$ Hence $\ell_0=0.$
\medskip
\end{vproof}

\noindent
\textbf{Acknowledgements} \\
The author is grateful to Professor J.~I.~D\'{\i}az for some useful discussions about this paper.

\baselineskip .4cm

\addcontentsline{toc}{section}{References}
\def\cprime{$'$}

\end{document}